\documentclass[12pt]{article}
\usepackage[utf8]{inputenc}
\usepackage[left=18mm, right=18mm, top=20mm, bottom=20mm]{geometry}
\usepackage{amsmath,amsfonts,amssymb,amsthm}
\usepackage{hyperref}
\usepackage{graphicx}
\usepackage{epstopdf}
\usepackage{color}
\usepackage{hyperref}
\usepackage{float}
\usepackage{endnotes}
\usepackage{ragged2e}
\theoremstyle{definition}
\newtheorem{exmp}{Example}[section]

\providecommand{\keywords}[1]
{
	\small	
	\textbf{\textit{Keywords---}} #1
}
\providecommand{\MSC}[1]
{
	\small	
	\textbf{\textit{MSC Subject Classification---}} #1
}
\title{Arc length based WENO scheme for Hamilton-Jacobi Equations}
\author{Samala Rathan$^{1}$\thanks{\textit{Email:rathans.math@iipe.ac.in}} , Biswarup Biswas$^{2}$\thanks{\textit{Email:biswarupb7@gmail.com}} \\
	\small $^{1}$\textit{Faculty of Mathematics, Indian Institute of Petroleum \& Energy-Visakhapatnam, India-530003} \\
	\small $^{2}$\textit{Department of Mathematics, Indian Institute of Technology-Delhi, New Delhi, India-110016} \\
}
\date{} 
\begin{document}
	\maketitle

		\begin{abstract}
In this article, novel smoothness indicators are presented for calculating the nonlinear weights of weighted essentially non-oscillatory scheme to approximate the viscosity numerical solutions of Hamilton-Jacobi equations. These novel smoothness indicators are constructed from the derivatives of reconstructed polynomials over each sub-stencil. The constructed smoothness indicators measure the arc-length of the reconstructed polynomials so that the new nonlinear weights could get less absolute truncation error and gives a high-resolution numerical solution. Extensive numerical tests are conducted and presented to show the performance capability and the numerical accuracy of the proposed scheme with the comparison to the classical WENO scheme.
		\end{abstract}
		
\\
\MSC{65M06, 65M12, 65M70, 41A10.}
		
\section{Introduction}
\label{S:1}
This article is concerned about the design of a new fifth-order weighted essentially non-oscillatory scheme to approximate the Hamilton-Jacobi (HJ) equations 
\begin{equation}\label{HJ}
\begin{aligned}
& \phi_t(\textbf{x},t)+H(\textbf{x}, t, \phi, D\phi)=0,\\
& \phi(\textbf{x},0)=\phi_0(\textbf{x}),
\end{aligned}
\end{equation}
where $\textbf{x}\in\mathbb{R}^{d},\,t>0$. Here, $H$ is known as the Hamiltonian and $D\phi=\left(\phi_{x},\phi_{y}\right)$. The physical significance of equation (\ref{HJ}) is important as it appears in several applications such as seismic waves, image processing, optimal control, calculus of variations, robotic navigation, crystal growth, etching, differential games and geometric optics. 
The difficulty in dealing with these equations (\ref{HJ}) is that it develops the discontinuous derivatives even with smooth initial data. As a result, the solutions for equation (\ref{HJ}) are not available as a unique sense and to further study solutions are understood in a weaker sense. Weak solutions are introduced by the notion of viscosity solutions\cite{CL, CEL, CIL}.
\par 
Indifference to the theoretical results, much attention gained to approximate the solutions of (\ref{HJ}) numerically. The authors, Crandall and Lions \cite{CL1}, first introduced the numerical approximation of first-order accurate monotone finite difference schemes to (\ref{HJ}) and its convergence theory is studied in \cite{PES}. An interesting fact that these equations (\ref{HJ}) are having the close connection to the conservation laws in its solution structure, so from observation, it is readily possible to obtain the exact solutions of these equations from those of conservation laws. So, an immediate consequence of the successful numerical methods for conservation laws can be adapted intently for solving (\ref{HJ}). In such direction, the authors Osher and Sethian \cite{OS} constructed a second-order essentially non-oscillatory (ENO) scheme, and Osher and Shu \cite{OSC} presented higher-order ENO schemes.  In 2000, Jiang and Peng \cite{JP} proposed a fifth-order finite difference weighted ENO, referred it as WENO-JP scheme, for solving HJ-equations. Further, the WENO schemes comes in to the literature to solve these equations, and many higher-order WENO variants developed, such as, Hermite WENO schemes on structured meshes by Qiu \cite{Qiu, Qiu1}, Qiu and Shu \cite{Qiu2} and Zheng and Qiu \cite{ZS1}; on unstructured meshes Lafon and Osher \cite{LO}, Abgrall \cite{AB}, Augoula and Abgrall \cite{AA}, Zhang and Shu \cite{ZS}, Li and Chan \cite{LC}. Few more notable WENO schemes are central WENO schemes \cite{BL1}, weighted power ENO schemes \cite{SQ}, mapped WENO schemes \cite{BL2, HLY}, symmetrical WENO schemes \cite{AAD} and WENO-ZQ scheme\cite{ZQ}. There are also central high-resolution schemes by Kurgunov and Tadmor \cite{KT}, Lin and Tadmor \cite{LT}, Bryson and Levy \cite{BL}. The discontinuous Galerkin (DG) framework also applied to \ref{HJ} by Hu and Shu \cite{HS}, Guo et al. \cite{GLQ}, Yan and Osher \cite{YO}, Cheng and Wang \cite{CW1} and Cheng and Shu \cite{CS2}.
\par
This work aims to design a new fifth-order WENO scheme to solve Hamilton-Jacobi equations. For the design of this new WENO scheme, we develop a set of new nonlinear weights. These nonlinear weights measure the smoothness of locally reconstructed polynomials over each sub-stencils. The design of the smoothness measurements is our main focus, and it is done by using the measuring the length of the curve of reconstructed polynomials over each sub-stencil. Such developments can be seen in \cite{BD} for the finite volume formulation for conservation laws; here, we introduced such concept in finite difference approach and extended to solve (\ref{HJ}). These smoothness indicators approximate the derivatives involved in the reconstructed polynomials effectively with the curve length based measurements. Such measurements lead to give the nonlinear weights which could get less truncation error and higher-resolution numerical solution near steep gradient regions, discontinuities and at complex structures. The resulted numerical scheme achieves the same fifth-order of accuracy as of classical WENO-JP scheme \cite{JP} but with less absolute errors. Numerical experiments are conducted in one and two-dimensions to show the effectiveness of the proposed numerical scheme.
\par
This paper is organized as follows. In Section 2, we start with a presentation of one-dimensional finite difference WENO reconstruction procedure for solving Hamilton-Jacobi equations, and then we have presented the proposed scheme in an algorithmic approach with the inclusion of time discretization scheme. Section 3 detailed out about the extensive numerical experiments for one and two-dimensions to verify the accuracy and robustness of the numerical scheme. Finally, concluding remarks are given in Section 4.  
\section{WENO scheme for HJ-Equations}
\label{S:2}
In this section, we will present a general framework of a fifth-order WENO scheme for solving Hamilton-Jacobi equations
\begin{equation}\label{HJ1}
\begin{aligned}
& \phi_t(x,t)+H(\phi_x(x,t))=0, \,\, (x,t) \in \Omega \times(0,\infty),\\
& \phi(x,0)=\phi_0(x),
\end{aligned}
\end{equation}
satisfying with suitable boundary conditions. We divide the computational domain $\Omega$ into uniformly distributed intervals $I_{i}$ with the mesh width $\Delta x=(x_{i-1}, x_{i})$ and with the mesh grid $x_{i}=i\Delta x$. We use the following notation
\begin{equation}
\phi_{i}=\phi(x_{i}), \,\, \Delta_x^{+}\phi_{i}=\dfrac{\phi_{i+1}-\phi_{i}}{\Delta x}, \,\, \Delta_x^{-}\phi_{i}=\dfrac{\phi_{i}-\phi_{i-1}}{\Delta x}.
\end{equation} 
Observe that $\phi_{i}^{n}$ and $\phi_{x,i}^{n}$ are the approximate values of $\phi(x_{i}, t^{n})$ and $\phi_x(x_{i}, t^{n})$ respectively, where $\phi_x$ indicate the derivative of $\phi$ with respect to $x.$ A semi-discrete conservative form for equation \eqref{HJ1} is
\begin{equation}\label{SDE}
\dfrac{d\phi_i(t)}{dt}=\mbox{{L}}(\phi_{i})=-\hat{H}\left(x_{i}, t, \phi_{i}, \phi_{x,i}^{+}, \phi_{x,i}^{-}\right),
\end{equation}
where $\phi_{x,i}^{+}$ and $\phi_{x,i}^{-}$ are the numerical approximation to $\phi_{x}(x_i)$. The Hamiltonian $\hat{H}$ is called as a numerical Hamiltonian and is a Lipschitz continuous monotone flux consistent with physical flux $H$, that is,
\begin{equation*}
\hat{H}(x,t,\phi,\phi_x,\phi_x)=H(x,t,\phi,\phi_x).
\end{equation*} Monotonicity refers that $\hat{H}$ is non-increasing in its fourth argument and non-decreasing in its fifth argument. We apply the Lax-Friedrichs flux 
\begin{equation}
\hat{H}(x,t,\phi,u^{+},u^{-})=H\left(x,t,\phi,\dfrac{u^{+}+u^{-}}{2}\right)-\alpha\left(\dfrac{u^{+}-u^{-}}{2}\right),
\end{equation}
where $u=\phi_x$ and $\alpha=\max_{u}|H_{1}(u)|$. The term $H_{1}$ represents the partial derivative of $H$ with respect to $\phi_x$. Now, we describe the reconstruction procedure for the approximation of $\phi_x$ from the left and right sides of the point $x_i$, that is, $\phi_{x,i}^{-}, \phi_{x,i}^{+}$ in order to calculate the approximate Hamiltonian $\hat{H}$ in an algorithm manner. 
\newline
\textbf{\textit{Step 1:}} To approximate $\phi_{x,i}^{-}$, we consider a fourth degree polynomial defined by $p_1^{-}(x)$ on larger stencil $S=\{x_{i-3},...,x_{i+2}\}$  which satisfies
\begin{equation}
\frac{1}{\Delta x}\int_{x_{j}}^{x_{j+1}}p_1^{-}(\xi)d\xi=\Delta_x^{+}\phi_{j}, j=i-3,...,i+1.
\end{equation}
The explicit expressions of fourth degree reconstructed polynomial $p_1^{-}(x_i)$ is
\begin{equation}
\begin{aligned}
\phi_{x,i}^{-}=p_1^{-}(x_i)= & \frac{1}{30}\Delta_x^{+}\phi_{i-3}-\frac{13}{60}\Delta_x^{+}\phi_{i-2}+\frac{47}{60}\Delta_x^{+}\phi_{i-1}-\frac{9}{20}\Delta_x^{+}\phi_{i}-\frac{1}{20}\Delta_x^{+}\phi_{i+1},
\end{aligned}
\end{equation}
which is a fifth-order approximation to $\phi_{x}({x_i})$. To approximate $\phi_{x,i}^{+}$, again we consider a fourth degree polynomial defined by $p_1^{+}(x)$ on larger stencil $T=\{x_{i-2},...,x_{i+3}\}$  which satisfies
\begin{equation}
\frac{1}{\Delta x}\int_{x_{j}}^{x_{j+1}}p_1^{+}(\xi)d\xi=\Delta_x^{+}\phi_{j}, j=i-2,...,i+2.
\end{equation}
The explicit expressions of fourth degree reconstructed polynomial $p_1^{+}(x_i)$ is
\begin{equation}
\begin{aligned}
\phi_{x,i}^{-}=p_1^{-}(x_i)= & \frac{-1}{20}\Delta_x^{+}\phi_{i-2}+\frac{9}{20}\Delta_x^{+}\phi_{i-1}+\frac{47}{60}\Delta_x^{+}\phi_{i}-\frac{13}{60}\Delta_x^{+}\phi_{i+1}-\frac{1}{30}\Delta_x^{+}\phi_{i+2},
\end{aligned}
\end{equation}
which is again a fifth-order approximation to $\phi_{x}({x_i})$. 
\newline
\textbf{\textit{Step 2:}} To construct a fifth-order WENO approximation to $\phi_{x,i}^{-}$ and $\phi_{x,i}^{+}$, we consider second degree polynomials $p_{k,1}^{-}(x)$ over each of the substencil $S^{-}_{k}=\{x_{i-3+k},...,x_{i+k}\}, \,\, k=0,1,2$ of stencil $S$ and $p_{k,1}^{+}(x)$ over each of the substencil $T^{+}_{k}=\{x_{i-2+k},...,x_{i+k}\},\,\,k=0,1,2$ of stencil $T$ respectively that satisfies
\begin{equation}
\begin{aligned}
& \frac{1}{\Delta x}\int_{x_{j}}^{x_{j+1}}p_{0,1}^{-}(\xi)d\xi=\Delta_x^{+}\phi_{j}, j=i-3,i-2,i-1, \\
& \frac{1}{\Delta x}\int_{x_{j}}^{x_{j+1}}p_{1,1}^{-}(\xi)d\xi=\Delta_x^{+}\phi_{j}, j=i-2,i-1,i, \\
& \frac{1}{\Delta x}\int_{x_{j}}^{x_{j+1}}p_{2,1}^{-}(\xi)d\xi=\Delta_x^{+}\phi_{j}, j=i-1,i,i+1. \\
\end{aligned}
\end{equation}
and
\begin{equation}
\begin{aligned}
& \frac{1}{\Delta x}\int_{x_{j}}^{x_{j+1}}p_{0,1}^{+}(\xi)d\xi=\Delta_x^{+}\phi_{j}, j=i-2,i-1,i, \\
& \frac{1}{\Delta x}\int_{x_{j}}^{x_{j+1}}p_{1,1}^{+}(\xi)d\xi=\Delta_x^{+}\phi_{j}, j=i-1,i,i+1, \\
& \frac{1}{\Delta x}\int_{x_{j}}^{x_{j+1}}p_{2,1}^{+}(\xi)d\xi=\Delta_x^{+}\phi_{j}, j=i,i+1,i+2. \\
\end{aligned}
\end{equation}
The explicit expressions of second degree reconstructed polynomial $p_{k,1}^{-}(x)$ for $k=0,1,2$ are
\begin{equation}
\begin{aligned}
\phi_{x,i}^{-,0}=p_{0,1}^{-}(x_i)= & \frac{1}{3}\Delta_x^{+}\phi_{i-3}-\frac{7}{6}\Delta_x^{+}\phi_{i-2}+\frac{11}{6}\Delta_x^{+}\phi_{i-1}, \\
\phi_{x,i}^{-,1}=p_{1,1}^{-}(x_i)= & -\frac{1}{6}\Delta_x^{+}\phi_{i-2}+\frac{5}{6}\Delta_x^{+}\phi_{i-1}+\frac{1}{3}\Delta_x^{+}\phi_{i},\\
\phi_{x,i}^{-,2}=p_{2,1}^{-}(x_i)= & \frac{1}{3}\Delta_x^{+}\phi_{i-1}+\frac{5}{6}\Delta_x^{+}\phi_{i}-\frac{1}{6}\Delta_x^{+}\phi_{i+1},
\end{aligned}
\end{equation}
where $p_{k,1}^{-}(x)$ is third-order approximations to $\phi_{x}({x_i})$. Note that $p_{k,1}^{+}(x)$ is obtained by the symmetric procedure with respect to the point $x_i$ for the reconstruction of $\phi_{x,i}^{-}$. 
\newline \textbf{\textit{Step 3:}} In smooth regions, a linear combination of 
$\phi_{x,i}^{-,0},\, \phi_{x,i}^{-,1}\, \phi_{x,i}^{-,2}$ with $d_0$, $d_1$ and $d_2,$ 
\begin{equation}
\begin{aligned}
\phi_{x,i}^{-}=d_0 \phi_{x,i}^{-,0}+d_1 \phi_{x,i}^{-,1}+d_2 \phi_{x,i}^{-,2},
\end{aligned}
\end{equation}
we get the linear weights as
\begin{equation}
d_0=\frac{1}{10},\,d_1=\frac{6}{10},\,d_2=\frac{3}{10}. 
\end{equation}
\newline \textbf{\textit{Step 4:}} In nonsmooth regions, the WENO methodology consists a convex combination or weighted average of $\phi_{x,i}^{-,k},\, k=0,1,2,$ as
\begin{equation}
\begin{aligned}
\phi_{x,i}^{-}=\omega_0 \phi_{x,i}^{-,0}+\omega_1 \phi_{x,i}^{-,1}+\omega_2 \phi_{x,i}^{-,2},
\end{aligned}
\end{equation}
to sustain both accuracy and non-oscillatory behavior due to oscillations, where $\omega_k$, $k=0,1,2$ are known as nonlinear weights which are defined as
\begin{equation}\label{WLOC1}
\omega_k=\dfrac{\alpha_k}{\sum_{l=0}^{2}\alpha_l}, \, \alpha_k=\dfrac{d_k}{\left(\epsilon+\beta_{k}\right)^{2}},
\end{equation}
where $ 0<\epsilon<< 1$ is introduced to prevent the denominator becoming zero and is set to be $10^{-6}$. The term $\beta_k$ represents the measurement of the smoothness of a local solution defined as
\begin{equation}\label{newbeta}
\beta_k=\left(\mathcal{L}_{[x_{i-1},\,x_{i}]}(p^k_i)\right)^2, \,(k=0,1,2).
\end{equation}
where $\mathcal{L}_{[a,b]}$ is the length of a given polynomial in $[a,b]$ defined by\begin{equation}\label{length}
\mathcal{L}_{[a,\,b]}(P):=\int_{a}^{b}\sqrt{1+{P'(x)}^2} dx.
\end{equation}
Since $p^k_i(x)$ are polynomials of degree two they can be written as, \begin{equation}
p^k_i(x)=a_i^k+b_i^k x+c_i^k x^2.
\end{equation}
However, $a_i^k$'s are not required in the calculation of arc length \eqref{length} as it contains only the derivative of the polynomial. The terms $b_i^k$'s and $c_i^k$'s are explicitly given as
\begin{flalign*}
b_i^0&=\frac{-6 \Delta x \phi _{i-3}+24 \Delta x \phi _{i-2}-30 \Delta x \phi _{i-1}+12 \Delta x \phi _i+6 x_i \phi _{i-3}-18 x_i \phi _{i-2}+18 x_i \phi _{i-1}-6 x_i \phi _i}{6 \Delta x^3},&&\\
c_i^0&=\frac{-3 \phi _{i-3}+9 \phi _{i-2}-9 \phi _{i-1}+3 \phi _i}{6 \Delta x^3},&&
\end{flalign*}

\begin{flalign*}
b_i^1&=\frac{6 \Delta x \phi _{i-1}-12 \Delta x \phi _i+6 \Delta x \phi _{i+1}+6 x_i \phi _{i-2}-18 x_i \phi _{i-1}+18 x_i \phi _i-6 x_i \phi _{i+1}}{6 \Delta x^3},&&\\
c_i^1&=\frac{-3 \phi _{i-2}+9 \phi _{i-1}-9 \phi _i+3 \phi _{i+1}}{6 \Delta x^3},&&
\end{flalign*}

\begin{flalign*}
b_i^2&=\frac{6 \Delta x \phi _{i-1}-12 \Delta x \phi _i+6 \Delta x \phi _{i+1}+6 x_i \phi _{i-1}-18 x_i \phi _i+18 x_i \phi _{i+1}-6 x_i \phi _{i+2}}{6 \Delta x^3},&&\\
c_i^2&=\frac{-3 \phi _{i-1}+9 \phi _i-9 \phi _{i+1}+3 \phi _{i+2}}{6 \Delta x^3}.&&
\end{flalign*}
Let us denote the primitive function as \begin{equation*}
\mathcal{I}[b,c,z]:=\begin{cases}
\displaystyle \frac{(b+2 c z) \sqrt{(b+2 c z)^2+1}+\sinh ^{-1}(b+2 c z)}{4 c}, & \text{if $c\neq0$},\\
\sqrt{1+b^2}z, & \text{if $c=0$}.
\end{cases}
\end{equation*}
Then the explicit lengths are given by
\begin{equation*}
\mathcal{L}_{[x_{i-1},\,x_{i}]}(p_i^k)=\mathcal{I}[b_i^k,c_i^k,x_{i}]-\mathcal{I}[b_i^k,c_i^k,x_{i-1}],\,\, (k=0,1,2).
\end{equation*}
The Taylor series expansion of the nonliear weights shows
\begin{equation}
\begin{aligned}
\omega_0 &=& \displaystyle\frac{1}{10}+O(\Delta x^2)\\
\omega_1 &=& \displaystyle\frac{3}{5}+O(\Delta x^2)\\
\omega_2 &=& \displaystyle\frac{3}{10}+O(\Delta x^2)
\end{aligned}
\end{equation}
which is that the nonlinear weights conveges to the linear weights with second order of accuracy. 
\newline \textbf{\textit{Step 5:}}
With these new smoothness indicators, the resulted nonlinear weights in WENO algorithm gives the desired fifth-order convergence when it is combined with the third-order TVD-Runge-Kutta scheme
\begin{equation}
\frac{du}{dt}=\mbox{{L}}(\phi),
\end{equation}
for vector $\phi^n$ which is as follows 
\begin{align*}
\phi_{0}&=\phi^n,\\
\phi_{1}& =\phi_0+\Delta t \mbox{L}(\phi_0),\\
\phi_{2}& =\frac{3}{4}\phi_0+\frac{1}{4}(\phi_1+\Delta t \mbox{L}(\phi_{1})),\\
\phi_{3}& =\frac{1}{3}\phi_0+\frac{2}{3}(\phi_2+\Delta t \mbox{L}(\phi_{2})),\\
\phi^{n+1}&=\phi_{3}.
\end{align*}
 \textbf{Note:} The two-dimensional algorithm is developed by using the dimension-by-dimension fashion and the major difference between the implementation of numerical schemes are on the subroutines for computing the corresponding nonlinear weights for the proposed scheme and WENO-JP scheme \cite{JP}.
\section{Numerical results}
In this section, we test the one- and two-dimensional test cases with the proposed scheme WENO-$\mathcal{L}$ and presented compared numerical results with the classical WENO scheme WENO-JP. We use the CFL number to be as $0.6$ to evaluate the numerical solution.
\begin{exmp}\label{eg1}
	We consider the linear advection equation \begin{equation}\label{linear}
\phi_t+\phi_x=0, \,\,-1<x<1
\end{equation}
with the initial state $\phi(x,0)=-\cos(\pi x)$. Convergence rates of the schemes WENO-JP and WENO-$\mathcal{L}$ at time $t=2$ are tabulated in Table \ref{tbl:lin_adv} which shows the WENO-$\mathcal{L}$ scheme has achieved its desired order of accuracy and comparable to WENO-JP.\\
\begin{table}[htb!]
	\begin{tabular}{cc}
	\centering
	\begin{tabular}{ccccc}
		&&WENO-$\mathcal{L}$&&\\
		\hline N & $L^\infty$ error & Order &  $L^1$ error  & Order \\ 
		\hline 20 & 2.74e-03 & ... & 3.87e-03 & ... \\ 
		40 & 8.02e-05 & 5.09 & 1.06e-04 & 5.19 \\ 
		80 & 2.39e-06 & 5.07 & 3.06e-06 & 5.12 \\ 
		160 & 7.25e-08 & 5.04 & 9.18e-08 & 5.06 \\ 
		320 & 2.23e-09 & 5.02 & 2.81e-09 & 5.03 \\
		\hline 
	\end{tabular}
	&
	\begin{tabular}{ccccc}
		&&WENO-JP&&\\
		\hline N & $L^\infty$ error & Order &  $L^1$ error  & Order \\ 
		\hline 20 & 3.06e-03 & ... & 3.25e-03 & ... \\ 
		40 & 8.28e-05 & 5.21 & 9.89e-05 & 5.04 \\ 
		80 & 2.42e-06 & 5.10 & 3.03e-06 & 5.03 \\ 
		160 & 7.32e-08 & 5.05 & 9.30e-08 & 5.03 \\ 
		320 & 2.25e-09 & 5.02 & 2.87e-09 & 5.02 \\ 
		\hline 
	\end{tabular}
\end{tabular}
\caption{Accuracy table of Example \ref{eg1} at time $t=2$.}
\label{tbl:lin_adv}
\end{table}
\end{exmp}

\begin{exmp}\label{eg2}
	Next we produce computational results for the equation \begin{equation}
\phi_t-\cos(\phi_x+1)=0, \,\,-1<x<1
\end{equation}
with the initial $\phi(x,0)=-\cos(\pi x)$ at time $t=0.5/\pi^2$. Table \ref{cos1D} shows the convergence rate of both the schemes are very similar. 
\begin{table}[htb!]
	\begin{tabular}{cc}
	\centering
	\begin{tabular}{ccccc}
		&&WENO-$\mathcal{L}$&&\\
		\hline N & $L^\infty$ error & Order &  $L^1$ error  & Order \\ 
		\hline 20 & 5.40e-04 & ... & 3.04e-04 & ... \\ 
		40 & 4.98e-05 & 3.44 & 1.65e-05 & 4.21 \\ 
		80 & 4.68e-06 & 3.41 & 7.37e-07 & 4.48 \\ 
		160 & 2.31e-07 & 4.34 & 2.86e-08 & 4.69 \\ 
		320 & 9.21e-09 & 4.65 & 1.00e-09 & 4.84 \\ 
		\hline
	\end{tabular}
	&
	\begin{tabular}{ccccc}
		&&WENO-JP&&\\
		\hline N & $L^\infty$ error & Order &  $L^1$ error  & Order \\
		\hline 20 & 5.22e-04 & ... & 2.76e-04 & ... \\ 
		40 & 5.32e-05 & 3.29 & 1.56e-05 & 4.15 \\ 
		80 & 4.65e-06 & 3.52 & 7.21e-07 & 4.44 \\ 
		160 & 2.32e-07 & 4.32 & 2.88e-08 & 4.65 \\ 
		320 & 9.24e-09 & 4.65 & 1.01e-09 & 4.83 \\
		\hline
	\end{tabular}
\end{tabular}
\caption{Accuracy table of Example \ref{eg2} at time $t=0.5/\pi^2$.}
\label{cos1D}
\end{table}
\end{exmp}

\begin{exmp}\label{eg3}
	We consider the 2D Burgers equation \begin{equation}
\phi_t+\frac{(\phi_x+\phi_y+1)^2}{2}=0, \,\,-2<x,y<2
\end{equation}
with the initial condition $\phi(x,y,0)=-\cos(\pi (x+y)/2)$. Results at time $t=0.5/\pi^2$ are given in Table \ref{tab:2DBu_accr} which shows that the WENO-$\mathcal{L}$ scheme achieved the fifth order accuracy and the convergence rates are better than the WENO-JS schemes. 
\begin{table}[htb!]	
	\begin{tabular}{cc}
	\centering
	\begin{tabular}{ccccc}
		&&WENO-$\mathcal{L}$&&\\
		\hline N & $L^\infty$ error & Order &  $L^1$ error  & Order \\ 
		\hline 20$\times$20 & 2.50e-03 & ... & 4.20e-04 & ... \\ 
		40$\times$40 & 1.43e-04 & 4.13 & 1.58e-05 & 4.73 \\ 
		80$\times$80 & 6.34e-06 & 4.50 & 5.32e-07 & 4.90 \\ 
		160$\times$160 & 2.16e-07 & 4.87 & 1.66e-08 & 5.00 \\ 
		320$\times$320 & 6.84e-09 & 4.98 & 5.16e-10 & 5.01 \\  
		\hline
	\end{tabular}
	&
	\begin{tabular}{ccccc}
		&&WENO-JP&&\\
		\hline N & $L^\infty$ error & Order &  $L^1$ error  & Order \\
		\hline 20$\times$20 & 2.59e-03 & ... & 5.19e-04 & ... \\ 
		40$\times$40 & 1.51e-04 & 4.10 & 1.86e-05 & 4.80 \\ 
		80$\times$80 & 6.64e-06 & 4.51 & 6.75e-07 & 4.78 \\ 
		160$\times$160 & 2.24e-07 & 4.89 & 2.29e-08 & 4.88 \\ 
		320$\times$320 & 7.07e-09 & 4.99 & 7.64e-10 & 4.90 \\ 
		\hline
	\end{tabular}
\end{tabular}
\caption{Accuracy table of Example \ref{eg3} at time $t=0.5/\pi^2$.}
\label{tab:2DBu_accr}
\end{table}
\end{exmp}

\begin{exmp}\label{eg4}
	In this example we consider the following equation:
\begin{eqnarray}
\phi_t-\cos(\phi_x+\phi_y+1)=0, \,\,-2<x,y<2\\
\phi(x,y,0)=-\cos(\pi (x+y)/2).\nonumber
\end{eqnarray}
Computations are performed up to time $t=0.5/\pi^2$ and the corresponding accuracy table is given in Table \ref{tab:2DCos_accr}. We have similar convergence rate for both the schemes.
\begin{table}[htb!]	
	\begin{tabular}{cc}
	\centering
	\begin{tabular}{ccccc}
		&&WENO-$\mathcal{L}$&&\\
		\hline N & $L^\infty$ error & Order &  $L^1$ error  & Order \\ 
		\hline 20$\times$20 & 9.20e-04 & ... & 1.68e-04 & ... \\ 
		40$\times$40 & 1.45e-04 & 2.67 & 1.35e-05 & 3.63 \\ 
		80$\times$80 & 1.89e-05 & 2.94 & 8.36e-07 & 4.01 \\ 
		160$\times$160 & 1.12e-06 & 4.07 & 3.48e-08 & 4.59 \\ 
		320$\times$320 & 4.13e-08 & 4.77 & 1.20e-09 & 4.86 \\  
		\hline
	\end{tabular}
	&
	\begin{tabular}{ccccc}
		&&WENO-JP&&\\
		\hline N & $L^\infty$ error & Order &  $L^1$ error  & Order \\
		\hline 20$\times$20 & 9.48e-04 & ... & 1.75e-04 & ... \\ 
		40$\times$40 & 1.48e-04 & 2.68 & 1.37e-05 & 3.67 \\ 
		80$\times$80 & 1.95e-05 & 2.92 & 8.64e-07 & 3.99 \\ 
		160$\times$160 & 1.16e-06 & 4.07 & 3.61e-08 & 4.58 \\ 
		320$\times$320 & 4.23e-08 & 4.78 & 1.24e-09 & 4.86 \\ 
		\hline
	\end{tabular}
\end{tabular}
\caption{Accuracy table of Example \ref{eg4} at time $t=0.5/\pi^2$.}
\label{tab:2DCos_accr}
\end{table}
\end{exmp}

\begin{exmp}\label{eg5}
	We solve the linear equation \eqref{linear} with $\phi(x,0)=\phi_0(x-0.5)$ where $\phi_0(x)$ is given by,\begin{equation*}
\phi_0(x)=-\left(\frac{\sqrt{3}}{2}+\frac{9}{2}+\frac{2\pi}{3}\right)\left(x+1\right)+
\begin{cases}\displaystyle 
2 \cos(3\pi x^2/2)-\sqrt{3}		 & x<-1/3;\\
\displaystyle\frac{3}{2}+3\cos (2\pi x)    & -1/3\leq x<0;\\
\displaystyle\frac{5}{2}-3\cos (2\pi x)    & 0\leq x<1/3;\\
\displaystyle\frac{28+4\pi+\cos(3\pi x)}{3}+6\pi x(x-1)  & x\geq 1/3.
\end{cases}
\end{equation*}
Results at time $t=2$ and $t=8$ are plotted in Figure \ref{fig:Eg5}. The WENO-$\mathcal{L}$ schemes gives better results than the WENO-JP scheme near the sharp corners. 
\begin{figure}[htb!]
	\includegraphics[scale=0.65]{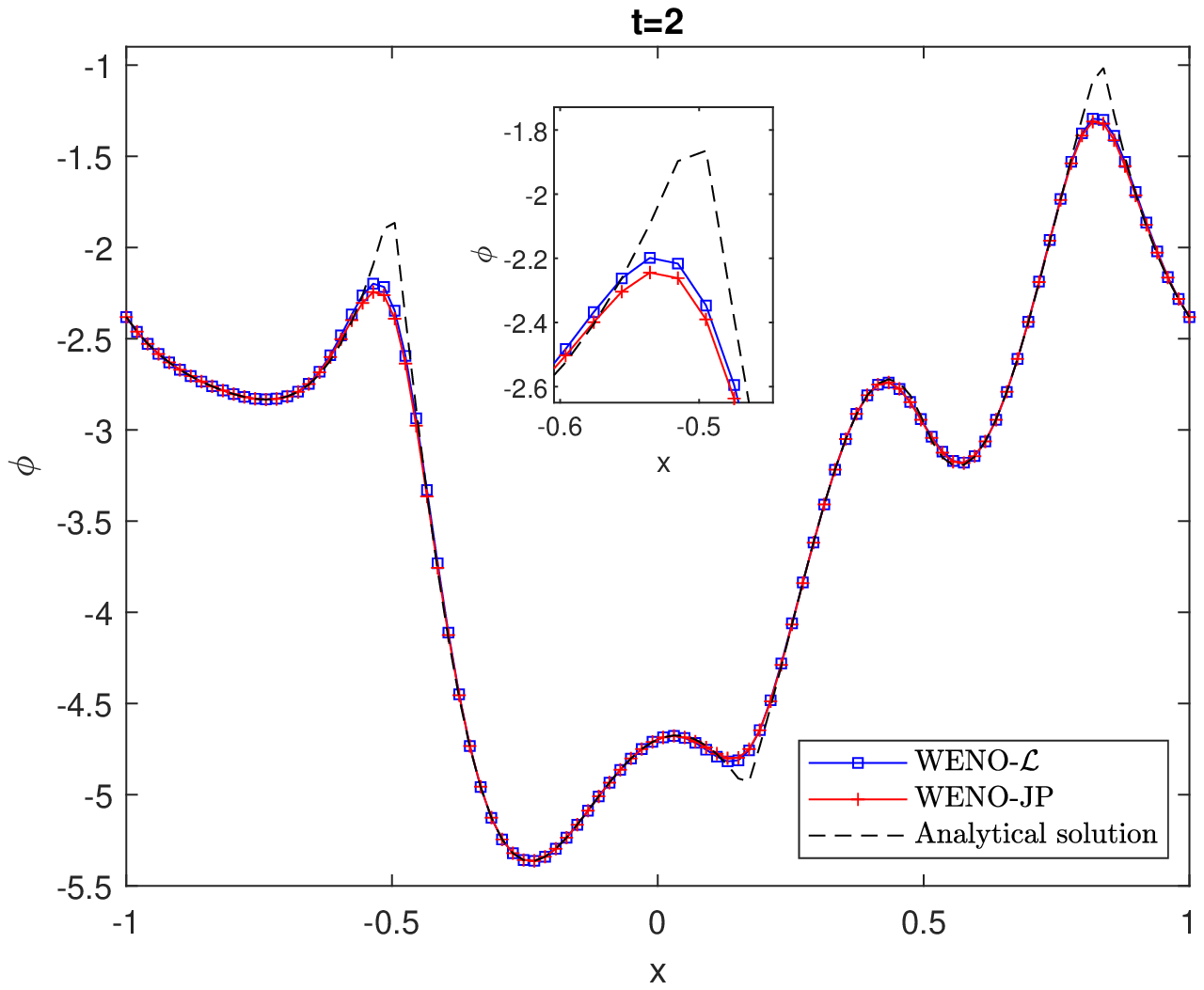}
	\includegraphics[scale=0.65]{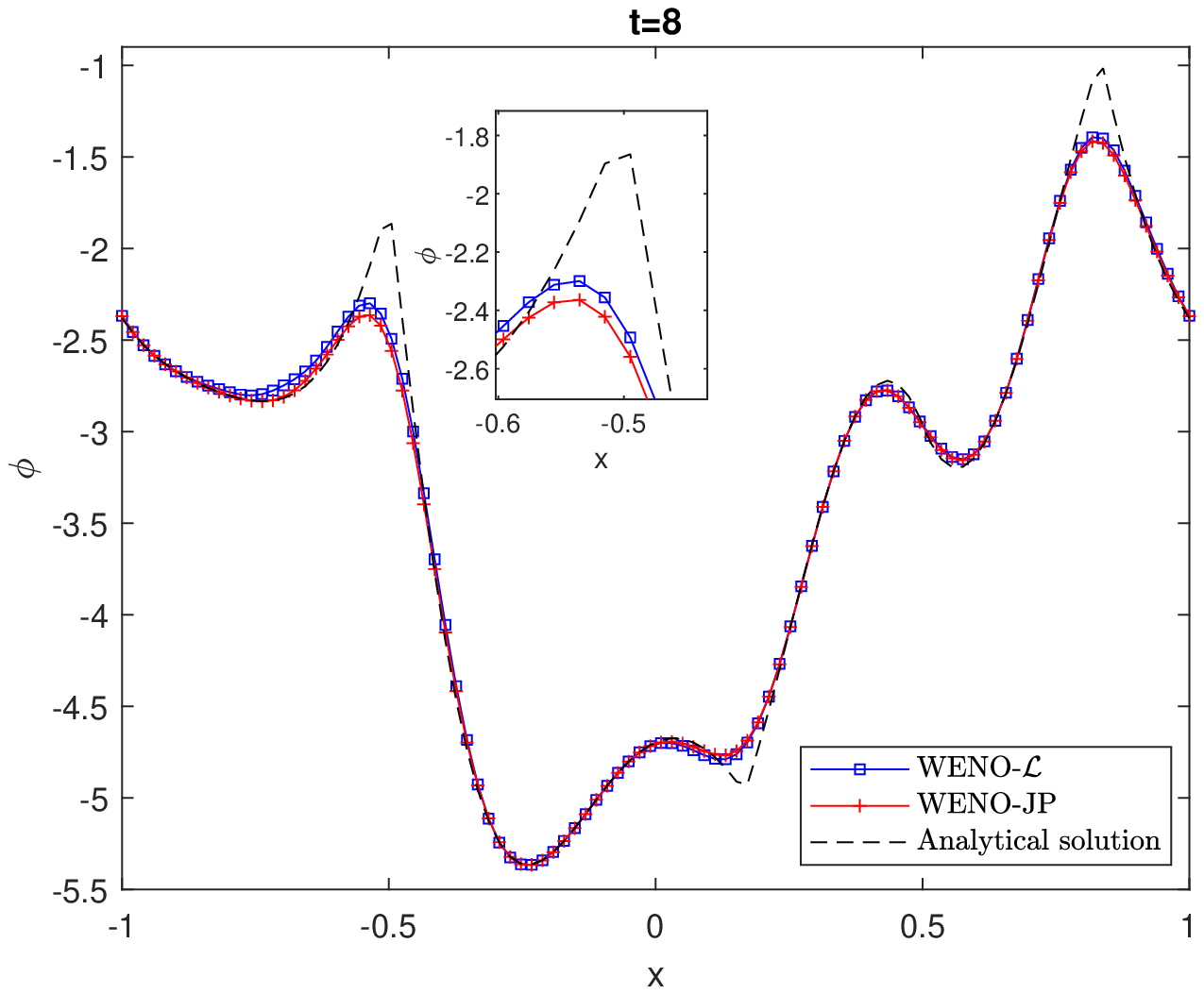}
	\caption{Solution of Example \ref{eg5} at time $t=2$(Left) and $t=8$(Right) using $N=100$.}
	\label{fig:Eg5}
\end{figure}
\end{exmp}

\begin{exmp}\label{eg6}
	We now consider the 1D Burgers equation \begin{equation}
\phi_t+\frac{(\phi_x+1)^2}{2}=0, \,\,-1<x<1
\end{equation}
with initial state $\phi(x,0)=-\cos(\pi x)$. Results at time $t=3.5/\pi^2$ are plotted in Figure \ref{fig:1DBurg} using 40 and 80 grids. Note that the discontinuous derivative is present in the solution at this time. Improvement on corners are clearly visible by the scheme WENO-$\mathcal{L}$. 
\begin{figure}[htb!]
	\includegraphics[scale=0.65]{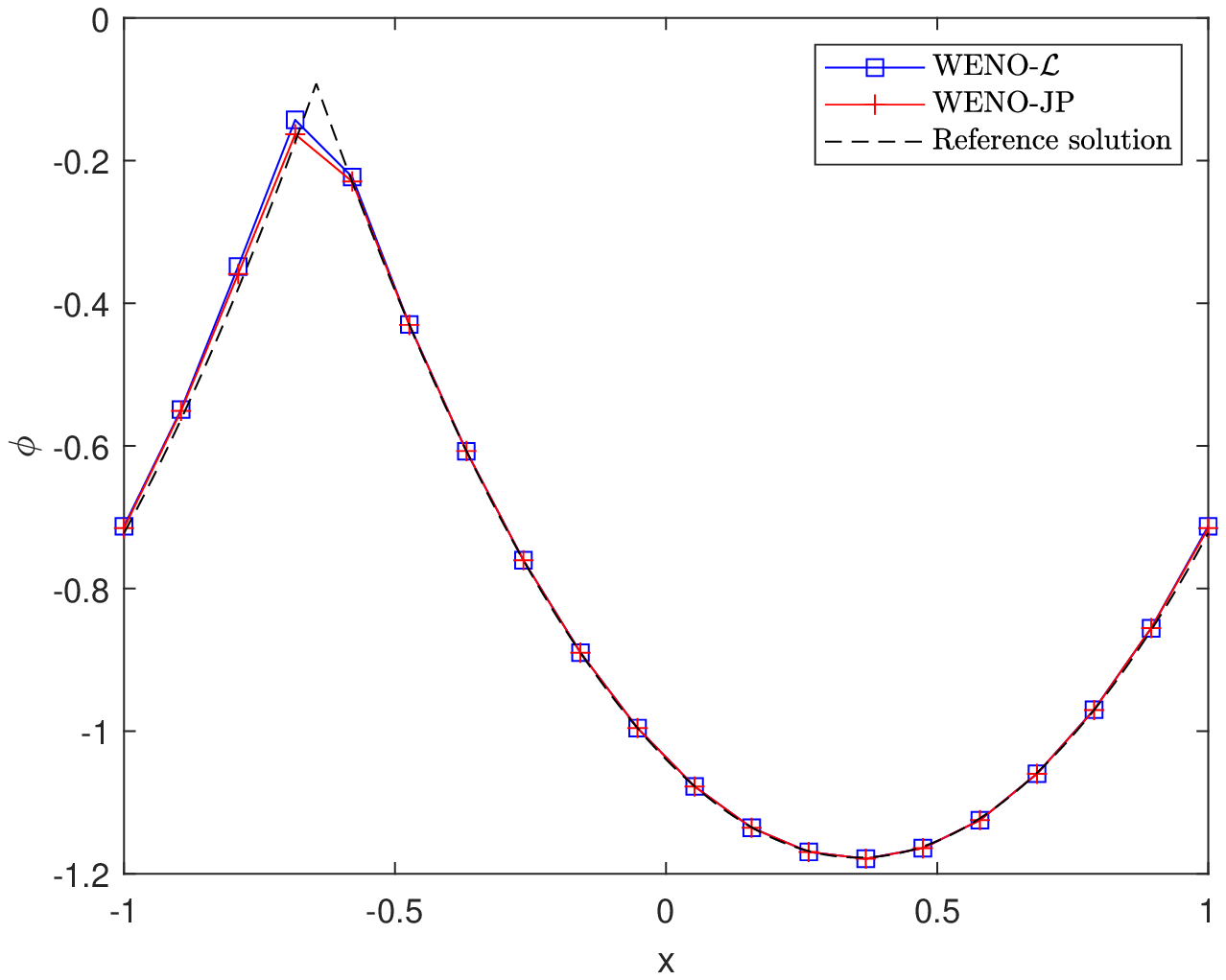}
	\includegraphics[scale=0.65]{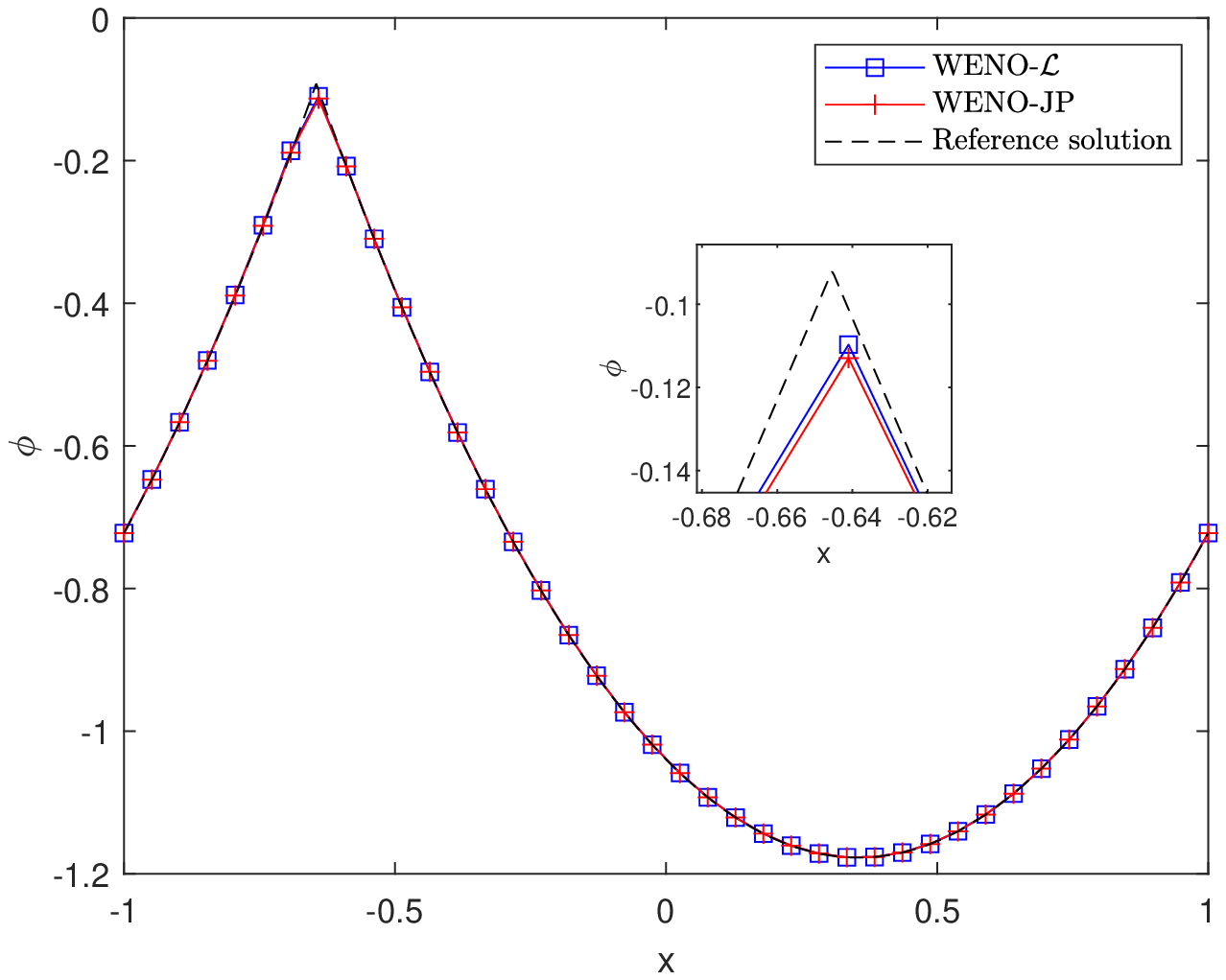}
	\caption{Solution of Example \ref{eg6} at time $t=3.5/\pi^2$ with grids $N=20$(Left) and $N=40$(Right).}	
	\label{fig:1DBurg}
\end{figure}
\end{exmp}

\begin{exmp}\label{eg7}
	We again consider the equation \eqref{cos1D}. This time we performed the computation up to time $t=1.5/\pi^2$. Results are displayed in Figure \ref{1dcos}. Both the schemes WENO-$\mathcal{L} $ and WENO-JP performed well in this case.
\begin{figure}[htb!]
	\includegraphics[scale=0.65]{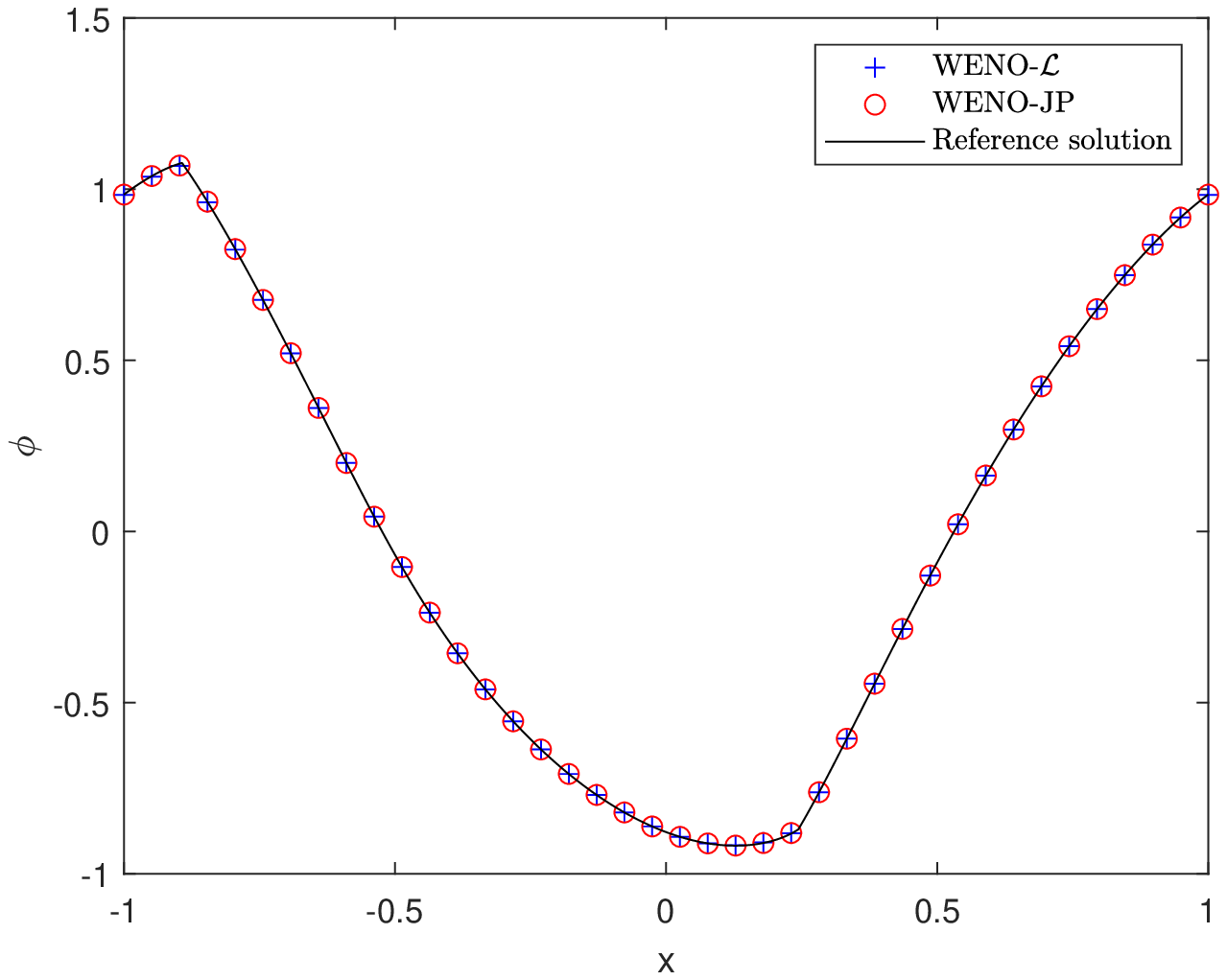}
	\includegraphics[scale=0.65]{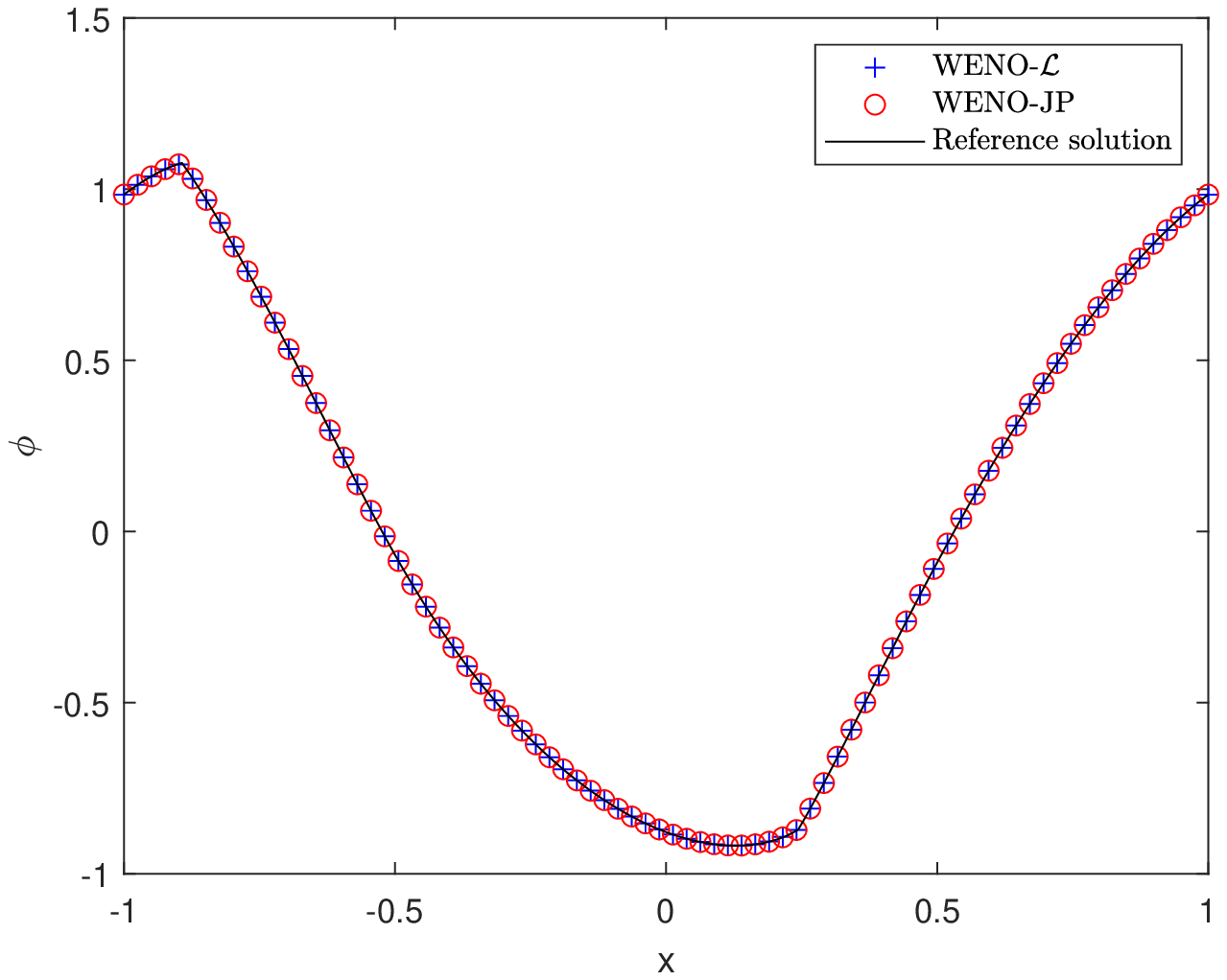}
		\caption{Solution of Example \ref{eg7} at time $t=3.5/\pi^2$ with grids $N=40$(Left) and $N=80$(Right).}
		\label{1dcos}
\end{figure}
\end{exmp}
\begin{exmp}\label{eg8}
	We solve the equation
	\begin{equation}
	\phi_t+\frac{1}{4}\left({\phi_x}^2-1\right)\left({\phi_x}^2-4\right)=0, \,\,-1<x<1
	\end{equation}
	with the initial condition $\phi(x,y,0)=-2|x|$. Results at $t=1$ are plotted in Figure \ref{fig:eg8} which shows that the top flatten geometry is captured well by both the schemes. 
	\begin{figure}[htb!]
			\centering
		\includegraphics[scale=0.75]{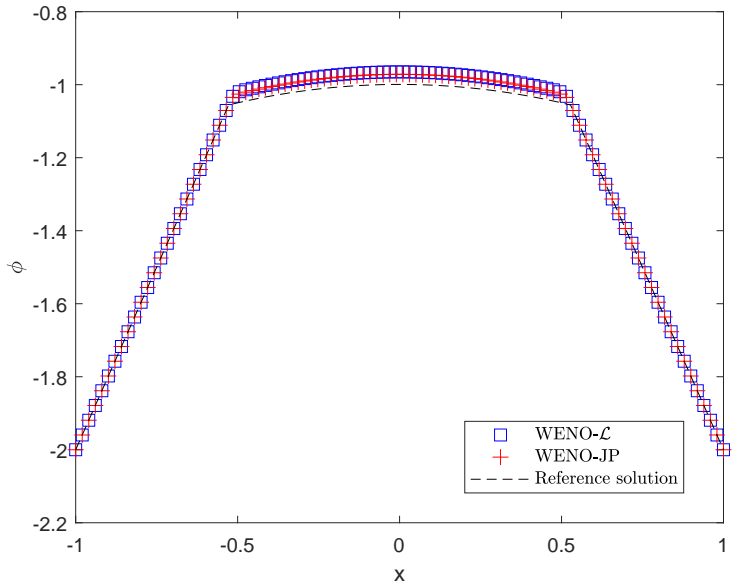}
		\caption{Solution of Example \ref{eg8} at time $t=1$ with grids $N=100$}
		\label{fig:eg8}
	\end{figure}
\end{exmp}

\begin{exmp}\label{eg9} 
	We consider the Example \eqref{eg3} and compute the solution at $t=1.5/\pi^2$. Results are obtained using $40\times 40$ mesh and displayed in Figure \eqref{fig:eg9}. Reasonably good resolution near sharp gradient is observed by the WENO-$\mathcal{L}$ scheme.	
	\begin{figure}[ht!]
		\includegraphics[scale=0.65]{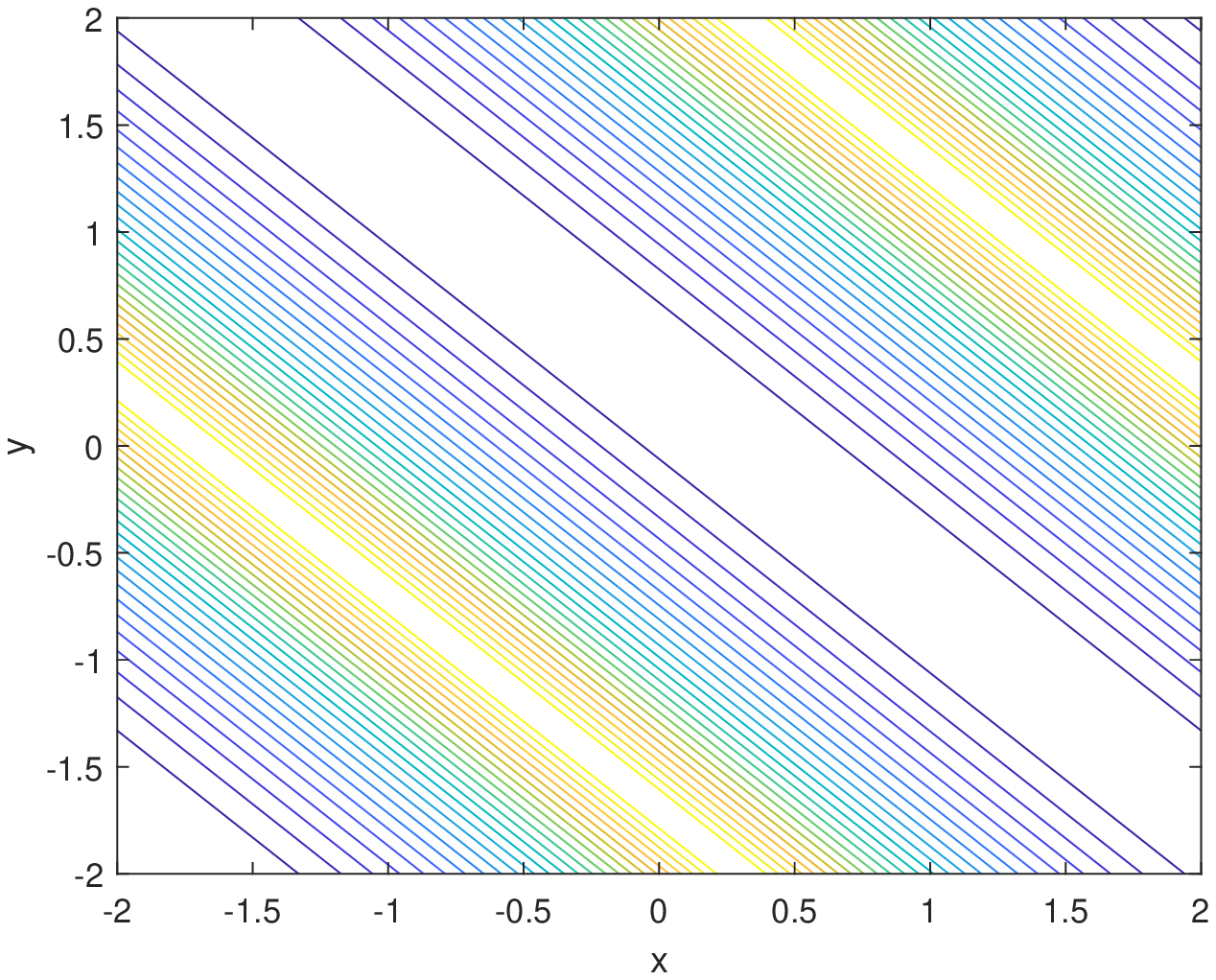}
		\includegraphics[scale=0.65]{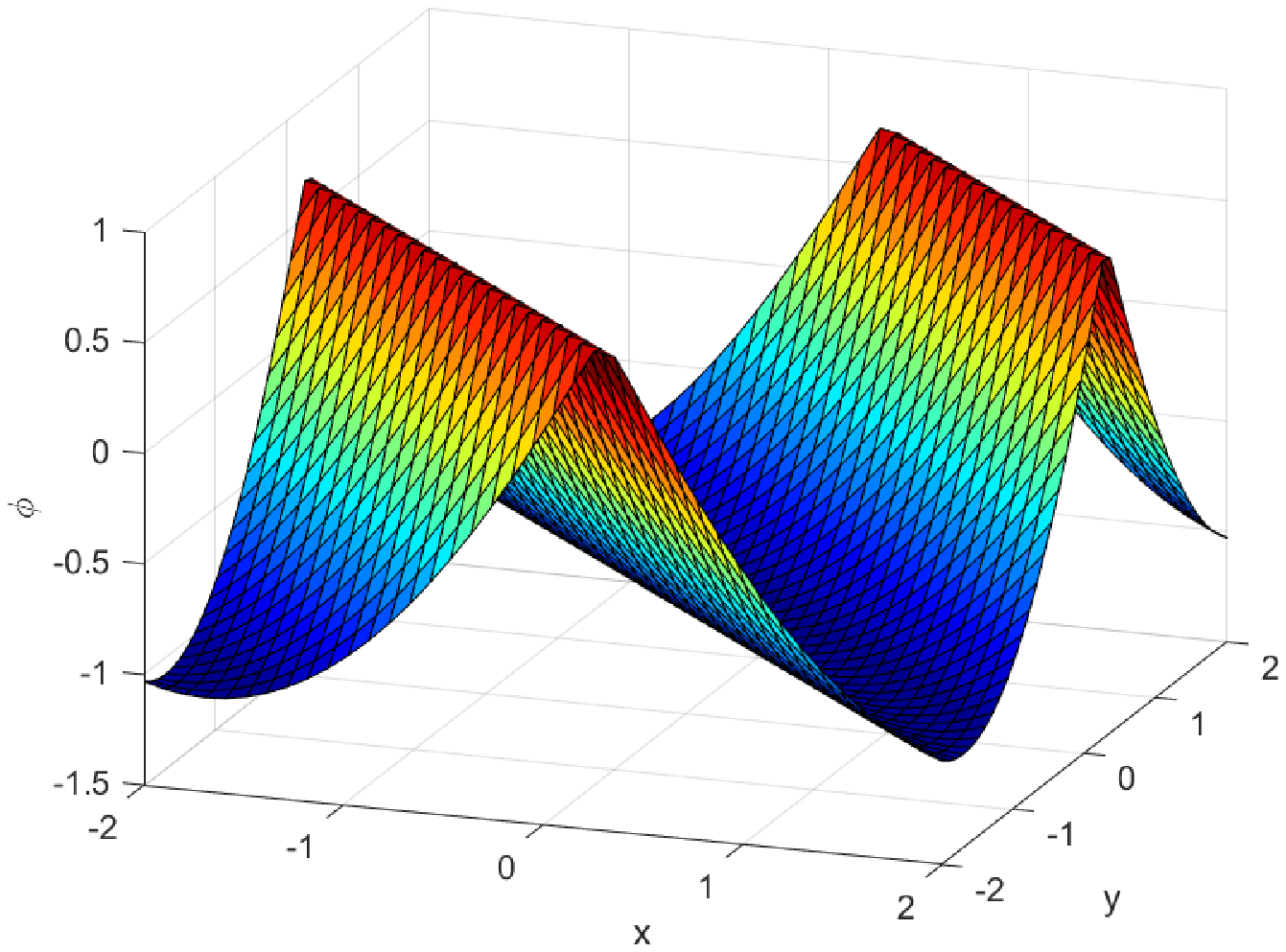}
		\caption{Contour (Left) and surface (Right) plot of the solution of Example \ref{eg9} at time $t=1.5/\pi^2$ using $40\times 40$ meshes.}
		\label{fig:eg9}
	\end{figure}
\end{exmp}
\begin{exmp}\label{eg10} 
	The following two dimensional equation with non-convex Hamiltonian is considered:
	\begin{equation}
\phi_t+\sin(\phi_x+\phi_y)=0,\,\,-1\leq x,y<1.
	\end{equation}
	The initial solution surface is $\phi(x,y,0)=\pi\left(|y|-|x|\right)$. Computational result is obtained at $t=1$ using $80\times80$ grids and Dirichlet boundary. We again observed sharp resolution in this case by WENO-$\mathcal{L}$ scheme.
	\begin{figure}[htb!]
		\includegraphics[scale=0.65]{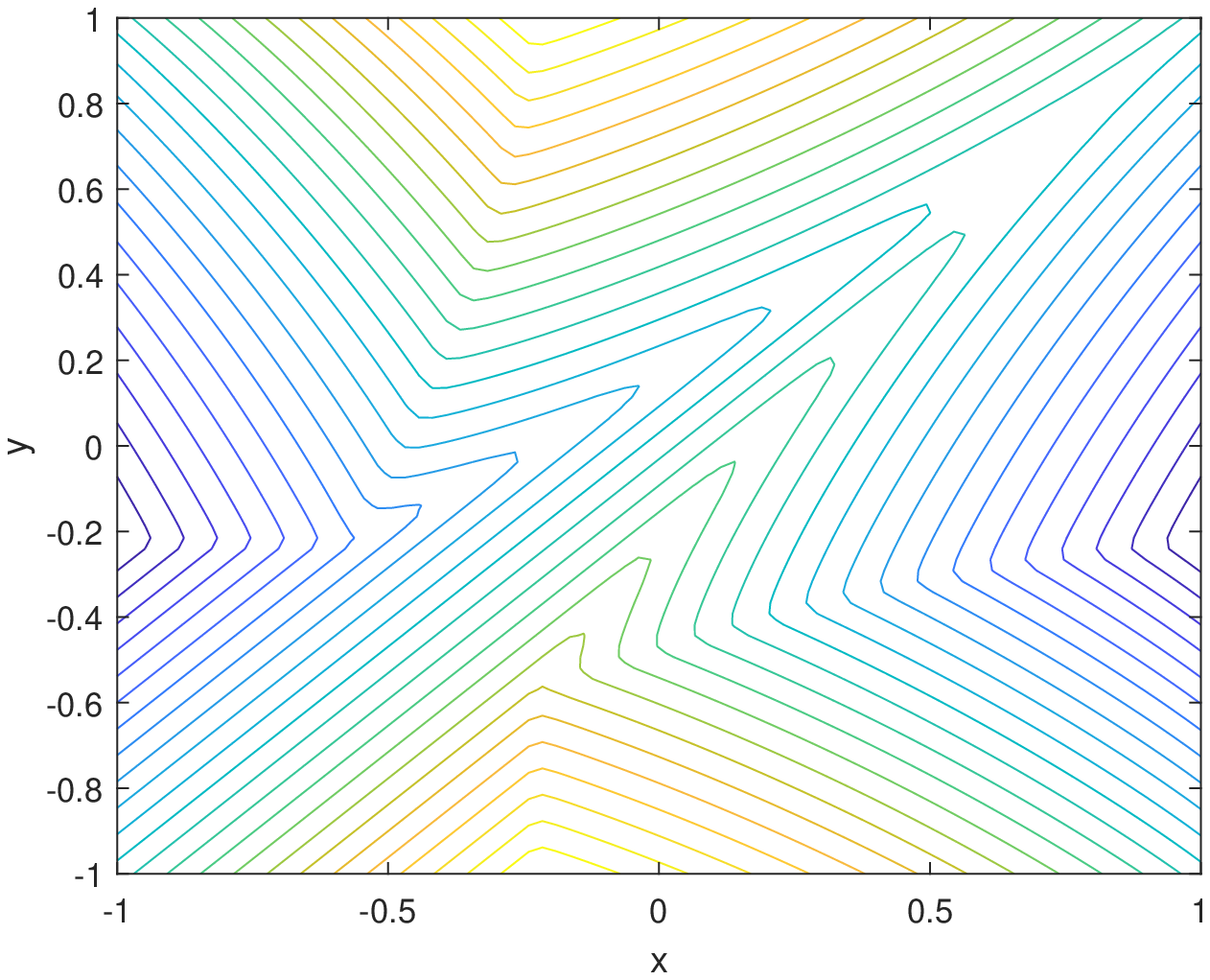}
		\includegraphics[scale=0.65]{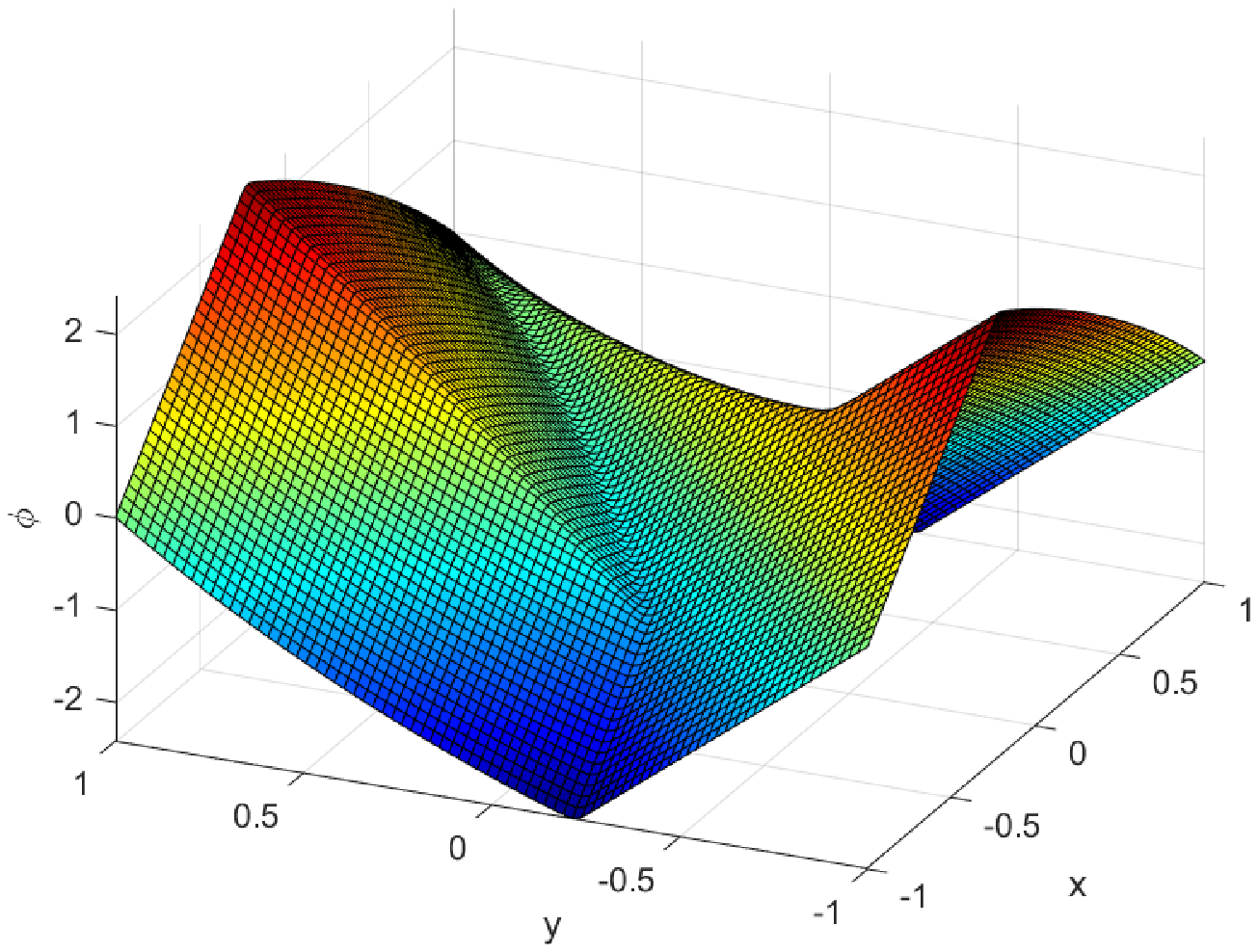}
		\caption{Contour (Left) and surface (Right) plot of the solution of Example \ref{eg10} at time $t=1$ using $80\times 80$ grid points.}
		\label{fig:eg10}
	\end{figure}
\end{exmp}
\begin{exmp}\label{eg11} 
We next consider a optimal control problem which is given by
\begin{equation}
\phi_t+\sin(y)\phi_x+\left(\sin(x)+sign\phi_y\right)\phi_y-\frac{1}{2}\sin^2(y)-(1-\cos(x))=0,\,\,-\pi\leq x,y<\pi,
\end{equation}
with $\phi(x,y,0)=0$. Solution is computed at $t=1$ and displayed in Figure \ref{fig:eg11}.
	\begin{figure}[htb!]
		\includegraphics[scale=0.65]{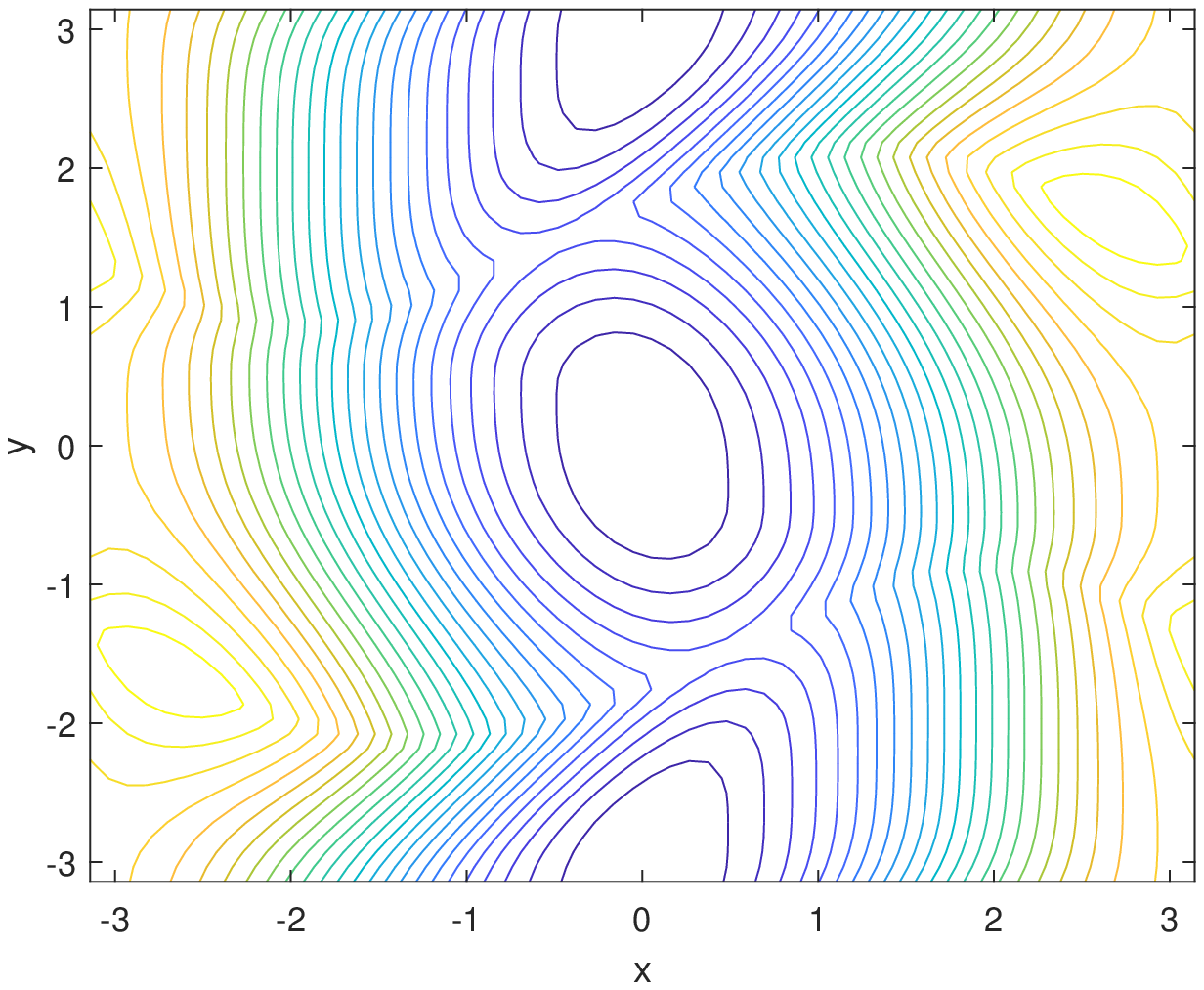}
		\includegraphics[scale=0.65]{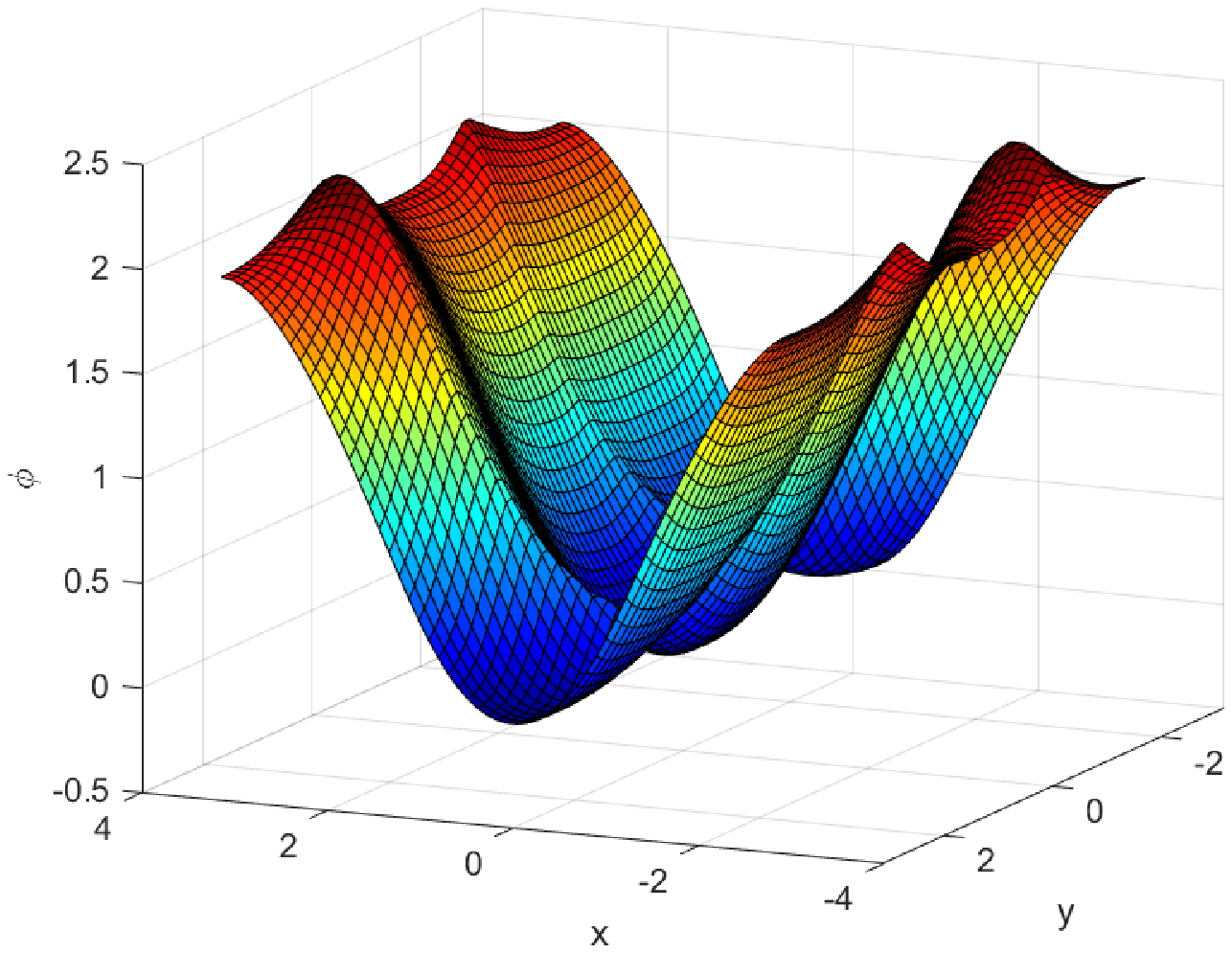}
		\caption{Contour (Left) and surface (Right) plot of the solution of Example \ref{eg11} at time $t=1$ using $60\times 60$ grid points.}
		\label{fig:eg11}
	\end{figure}
\end{exmp}
\begin{exmp}\label{eg12} 
	We solve the following 2D Eikonal equation which arises in geometric optics (see \cite{Jin1998}):
	\begin{eqnarray}
	\phi_t+\sqrt{{\phi_x}^2+{\phi_y}^2+1}=0,\,\,0\leq x,y<1;\\
	\phi(x,y,0)=\frac{1}{4}\left(\cos(2\pi x)-1\right)\left(\cos(2\pi y)-1\right)-1\nonumber
	\end{eqnarray}
Solution is computed at $t=0.6$ with $80\times 80$ grid points. Figure \ref{fig:eg12} shows very good resolution of the WENO-$\mathcal{L}$ scheme.
	\begin{figure}[htb!]
		\includegraphics[scale=0.65]{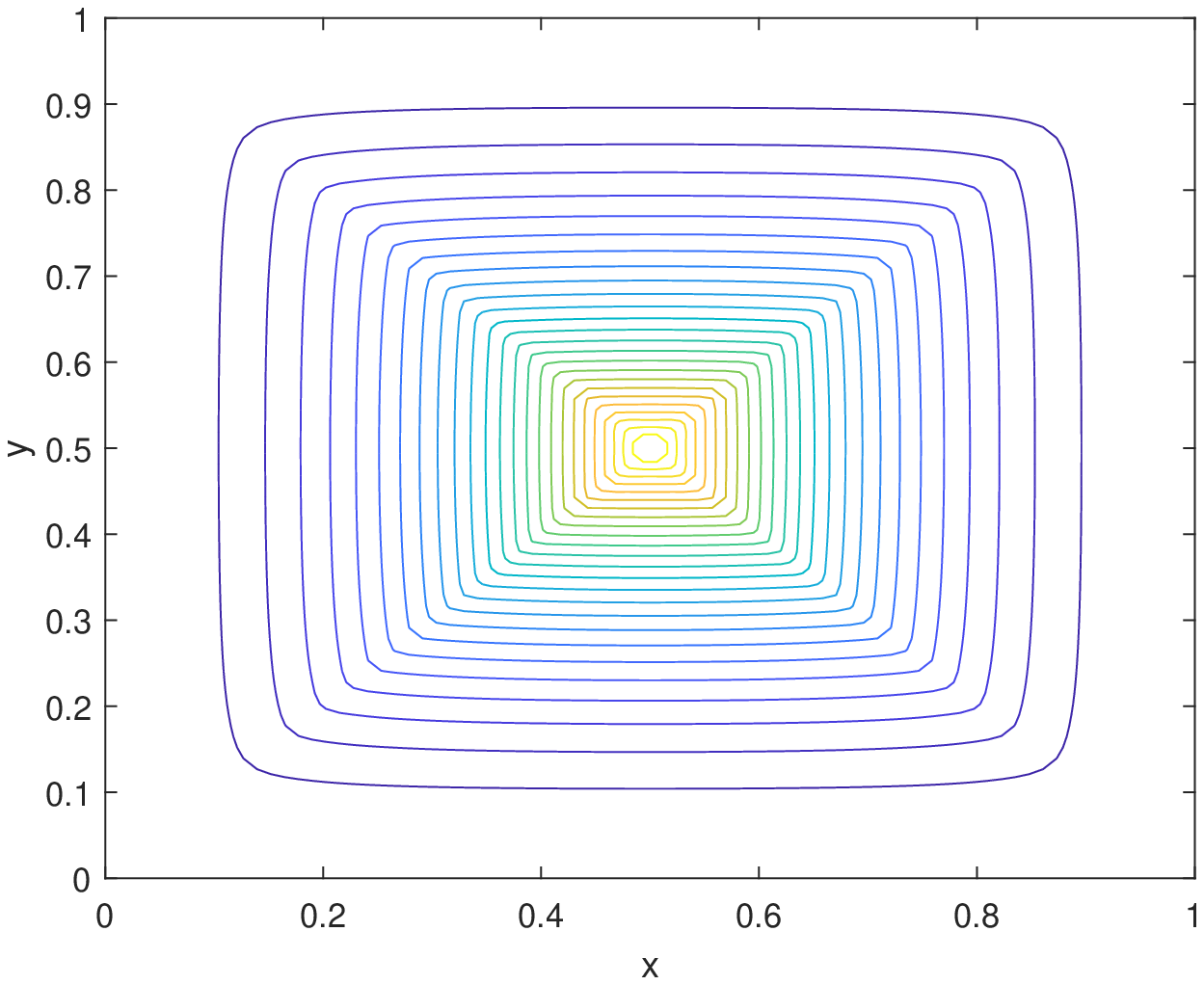}
		\includegraphics[scale=0.65]{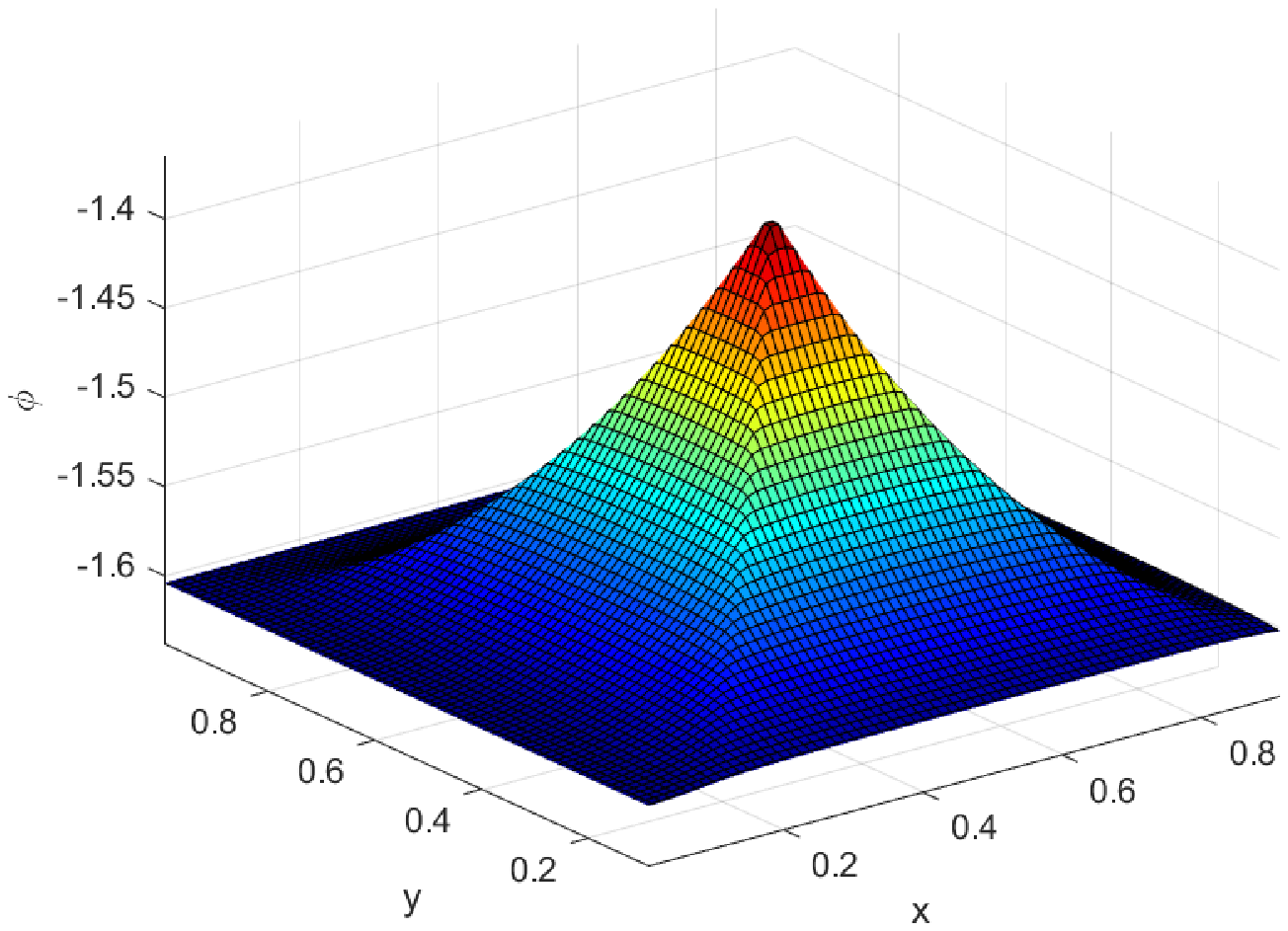}
		\caption{Contour (Left) and surface (Right) plot of the solution of Example \ref{eg12} at time $t=0.6$ using $80\times 80$ grid points.}
		\label{fig:eg12}
	\end{figure}
\end{exmp}

\begin{exmp}\label{eg13} 
	Now, we consider the following equation
	\begin{eqnarray}
	\phi_t-\left(1-\epsilon K\right)\sqrt{{\phi_x}^2+{\phi_y}^2+1}=0,\,\,0\leq x,y<1;\\
	\phi(x,y,0)=1-\frac{1}{4}\left(\cos(2\pi x)-1\right)\left(\cos(2\pi y)-1\right)\nonumber\\
	K=\frac{\phi_{xx}(1+\phi_{y})^2-2 \phi_{xy} \phi_{x} \phi_{y}+\phi_{yy}(1+\phi_{x})^2}{({\phi_x}^2+{\phi_y}^2+1)^{3/2}}\nonumber
	\end{eqnarray}
with $\epsilon=0$ which is a purely a convetion part and $\epsilon=0.1$ with the mesh grid points $60 \times 60$. It is obeserved that the WENO-$\mathcal{L}$ scheme produce a good resolution numerical solution in comparision to the WENO-JP scheme.
	\begin{figure}[htb!]
		\includegraphics[trim=40 0 90 10,clip,width=18cm,height=18cm]{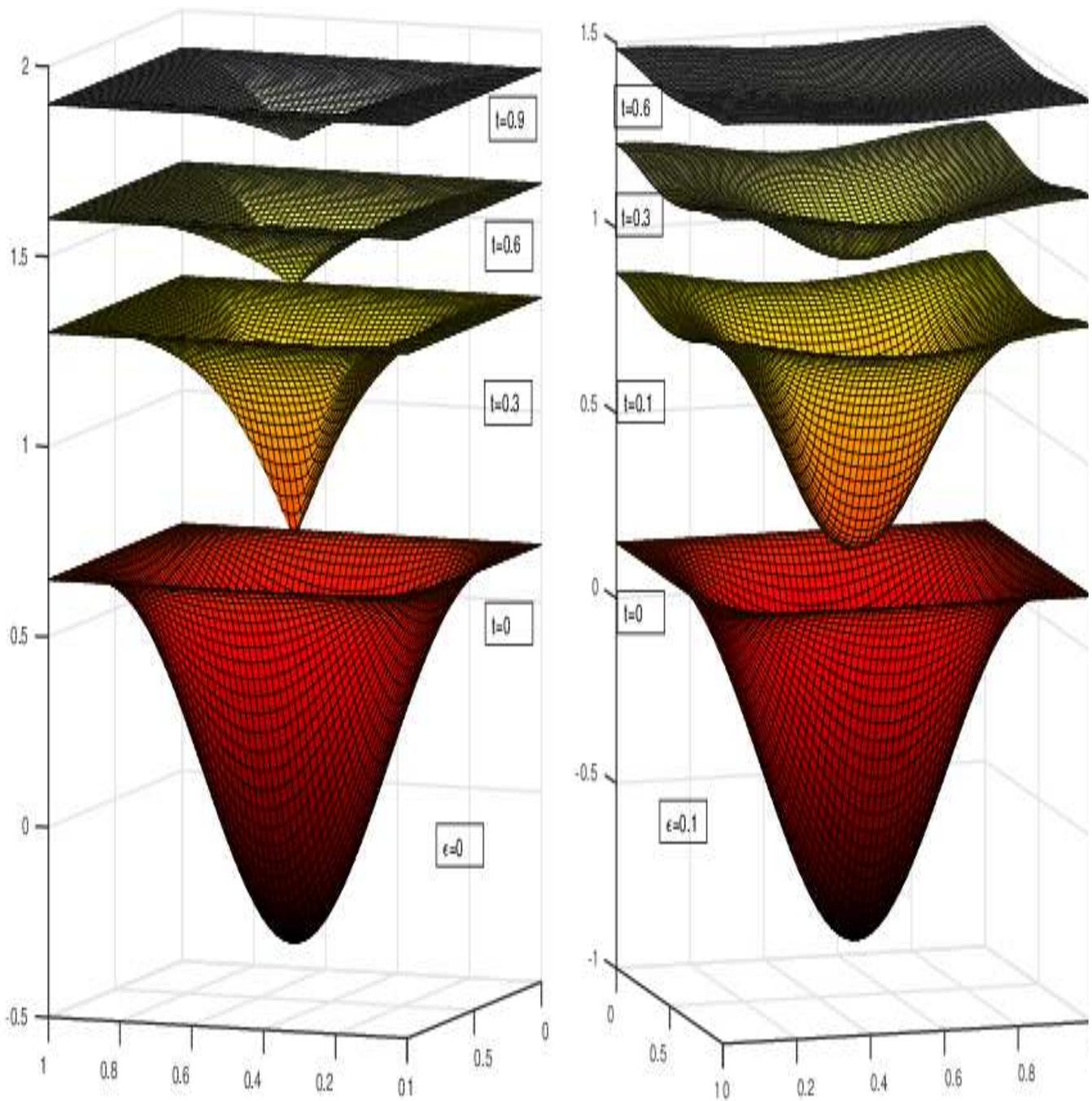}
		\caption{Surface plot left ($\epsilon=0$) and right ($\epsilon=0.1$) of the solution of Example \ref{eg13} using $60\times 60$ grid points.}
		\label{fig:eg13}
	\end{figure}
\end{exmp}
\section{Conclusions}
In this article, we have presented a new smoothness indicators in calculating the nonlinear weights for fifth-order weighted essentially non-oscillatory scheme to solve Hamilton-Jacobi equations. These new smoothness indicators are calculated based on derivatives of reconstructed polynomials over each sub-stencil that measures the length of the curve which fits on each cell. Extensive numerical tests are conducted to show the performance enhancement and the numerical accuracy of the proposed scheme. From the numerical results, it is observed that the proposed numerical scheme WENO-$\mathcal{L}$ produce higher resolution numerical solution in comparision to the classical WENO-JP scheme.

\end{document}